\documentclass[graybox]{svmult}

\usepackage{mathptmx}       
\usepackage{helvet}         
\usepackage{courier}        
\usepackage{type1cm}        
\usepackage{makeidx}         
\usepackage{graphicx}        
\usepackage{multicol}        
\usepackage[bottom]{footmisc}
\usepackage{amsmath}

\makeindex             

\begin{document}

\title*{On the occupation times of Brownian excursions and Brownian loops}
\author{Hao Wu}
\institute{D\'epartement de Math\'ematiques, Universit\'e Paris-Sud, 91405 Orsay Cedex, France\\
\email{hao.wu@math.u-psud.fr} \\
{The author acknowledges the support from a Fondation CFM-JP Aguilar grant.}}

%
%
\maketitle

\abstract{We study properties of occupation times by
 Brownian excursions and Brownian loops in two-dimensional domains.
This allows for instance to interpret some Gaussian fields, such as the Gaussian Free Fields as (properly
normalized) fluctuations of the total occupation time of a Poisson
cloud of Brownian excursions when the intensity of the cloud goes to
infinity.}

\newcommand{\mathbb}{}
\newcommand{\eps}{\epsilon}
\newcommand{\U}{\mathbb{U}}
\newcommand{\MR}{MR}

\bigbreak
\noindent
\textbf{Keywords:} Conformal invariance, Brownian excursion measure,
Brownian loop measure, Green's function.

\section{Introduction}

Conformal invariance of planar Brownian motion has been derived and exploited long ago by Paul L\'evy \cite {Levy}. See also B. Davis (Annals of Proba 1979) in particular his
derivation of Picard's big theorem. More recently, conformal invariance turned out to be an instrumental idea in the
study of various critical models from statistical physics in the
plane (see for instance \cite {Lawler2005, wernerleshouches} and the references therein).
Two basic important conformally invariant measures on random geometric objects
are the Brownian excursion measure and the
Brownian loop measure. Let us now very briefly describe these measures and the meaning of
conformal invariance relatively to these measures. For each open domain $D$ with non-polar
boundary in the plane, one can define these two measures in $D$
respectively denoted by $\mu_{D}$ and $\lambda_{D}$. These are
infinite but $\sigma$-finite measures on Brownian-type paths with
particular properties:
\begin{itemize}
\item $\mu_{D}$ is supported on the set of
Brownian excursions $(B_{t},t\leq\tau)$ in $D$ i.e. Brownian paths such
that $B_{0}$ and $B_{\tau}$ are in $\partial D$, while
$B(0,\tau)\subset D$.
\item $\lambda_{D}$ is supported on the set of Brownian loops
$(B_{t},t\leq\tau)$ i.e. Brownian paths in $D$ such that
$B_{0}=B_{\tau}$.
\end{itemize}
In fact, in both cases, it is useful to view these paths up to
monotone reparametrization (in the loop-case, one views the time-set modulo $\tau$ i.e., there is no ``starting point'' on the loop). Then, it turns out (see \cite{soup2004},\cite{conf2005} for details)
that for any conformal map $\Phi$ from $D$ onto $\Phi(D)$, the image
measures of $\mu_{D}$ and $\lambda_{D}$ under $\Phi$ are exactly
$\mu_{\Phi(D)}$ and $\lambda_{\Phi(D)}$.

These two measures on loops and on excursions
allow in some sense to get rid of the dependence of the measure on Brownian paths with respect to its starting point, see for instance the discussion in \cite {wernerleshouches}.

In the present text, we shall focus on the following type of results
(here and in the sequel, $dx$ or $dy$ will denote the area measure, and $x$ or $y$ will always denote points in the plane):
\begin{proposition}\label{main pro}
Suppose that $D$ is a simply connected domain and that $A$ and $B$ are two
open proper subsets of $D$. Then,
\begin{equation}\label{exc cov}
\mu_{D}(\int_{0}^{\tau}ds 1_{A}(\gamma_{s})\int_{0}^{\tau}ds
1_{B}(\gamma_{s}))=4\int_{A\times B}dx\ dy
\ G_{D}(x,y)
\end{equation}
and
\begin{equation}\label{loop cov}
\lambda_{D}(\int_{0}^{\tau}ds 1_{A}(\gamma_{s})\int_{0}^{\tau}ds
1_{B}(\gamma_{s}))=\int_{A\times B}dx\ dy
\ (G_{D}(x,y))^2,
\end{equation}
where $(\gamma_s, 0\leq s\leq \tau)$ is a Brownian excursion in $(\ref{exc cov})$ and a Brownian loop in $(\ref{loop cov})$, $G_D(x,y)$ denotes the usual Green's function in $D$ (with Dirichlet boundary conditions).
\end{proposition}
The Brownian excursion measure and the loop measure are infinite
measures, but they can be used to define random conformally
invariant collections of excursions and loops (i.e. under a
probability measure) by a Poissonization procedure. As explained in
\cite {wernerleshouches}, both these Poissonian clouds are of interest and useful in the context of
random planar conformally invariant curves of SLE-type: The ``excursion
clouds'' give rise to the restriction measures \cite {conf2005}, while the
``loop-soups (loop clouds)" are related to Conformal Loop Ensembles (see \cite {CLE}).

It is natural to study the cumulative occupation time of
these random collections of Brownian paths. The previous proposition can
then be viewed as a description of the covariance structure of these
cumulative occupation times (even if as we shall explain later, things are slightly more complicated in the
case of the loop measure because cumulative occupation times are infinite, so that a renormalization procedure is needed). By the classical central limit theorem,
in the asymptotic regime where the intensity of these
clouds goes to infinity, the fluctuations of these occupation times
converge (if properly normalized of course) to a Gaussian process with the
same covariance structure. This will in particular enable us to interpret the Gaussian
Free Field in terms of fluctuations of occupation times of high-intensity clouds of Brownian excursions.

Note that in \cite {LJ}, a different and more direct (as it involves no asymptotic) relation between the loop-soup occupation times and the
Gaussian Free Field (or rather its square) is pointed out.

Here is how the present paper is structured: In Section 2, we review various very elementary facts concerning Green's functions, their conformal invariance and their relation to Brownian motion and the Gaussian Free Field. In Section 3, we recall the definition of the Brownian excursion measure, we derive
 (\ref {exc cov}) and deduce from it the interpretation of the Gaussian Free Field as asymptotic fluctuations of the Excursions occupation time measure. In passing, we note a representation of the solution to the standard Dirichlet problem using Brownian excursions, that does not seem so well-known despite its simplicity.
Section 4 is the counterpart of Section 3 for Brownian loops instead of Brownian excursions. Finally, in Section 5, we
briefly mention a generalization of the previous results using some clouds of interacting pairs of excursions (via their intersection local-time) that
exhibits some relations between loops and excursions.

We will focus on two-dimensional domains, but many of our statements (in particular those on Brownian excursions) are also valid in higher dimensions. However, as the reader will see, we choose to base our proofs on conformal invariance, so that another approach would be needed to derive the results in dimensions greater than two.
We should also point out that the statements are in fact valid in non-simply connected domains, but again, some of our proofs,
in particular those dealing with the loop-measure,  would need to be changed in order to cover
non-simply connected planar domains (as we will use explicit expressions for the unit disc).

\medbreak
\noindent
\textbf{Acknowledgement:} This paper is based on my Master's thesis and was completed under the guidance
of my supervisor Professor Wendelin Werner.

\section{Review of basic notions}
\subsection{Generalities}
We first recall some classical facts about Brownian motion and its relation to harmonic functions, see for instance \cite {Durrett1984, Charles1978, Rao1977}
for further details or background.

Suppose that $D$ is a bounded planar domain, and that it has a smooth boundary.
Then, for any point $x$ in $D$, the distribution of the exit position from $D$ by a Brownian motion started at $x$ has a
continuous density with respect to the surface measure $\sigma (dz)$ on $\partial D$, called the {\em Poisson kernel}, that we will
denote by $h_D (x, z)$ for $z \in \partial D$. In other words, the exit distribution is $h_D (x,z) \sigma (dz)$.

This Poisson kernel is closely related to the solutions of the {\em Dirichlet problem} in $D$ (i.e.,
to find a harmonic function $u$ in $D$, that is continuous on $\overline {D}$
and equal to some prescribed continuous function $f$ on the boundary
of $D$). Indeed, the solution to the Dirichlet problem, if it exists, is given by
$$u(x) = \int\sigma (dz) h_D (x,z) f(z)  = E_x ( f( Z_\tau))$$
where $Z$ is a planar Brownian motion started from $x$ under the probability measure $P_x$ and $\tau$ denotes its exit time from $D$.

\medbreak

The {\em Green's function} in $D$, is the unique function in $D \times D$, such that for each $x \in D$,
$ y \mapsto G_D (x,y)$ is harmonic, vanishes on $\partial D$, and satisfies $G_D (x,y) \sim \pi^{-1} \log (1/|x-y|)$
when $y \to x$.

Alternatively, one can think of $G_{D}(x,y)dy$ as the expected time spent by $Z$
in the infinitesimal neighborhood of $y$ before exiting $D$. More precisely, if $A$
denotes an open set, the expected time spent by the Brownian motion $Z$ (started from $Z_0=x$) in $A$
before exiting $D$ is
$$
E_{x}(\int_{0}^{\tau}dt 1_{A}(Z_{t}))=\int_{A}dy
G_{D}(x,y).
$$

The Green's function is closely related to the  {\em Poisson problem} (i.e. to find a $C^{2}$ function $u$ in
$D$ such that $\Delta u=-2g$, where $g$ is some given continuous
function in $D$, with the property that $u$ is continuous on
$\overline{D}$ and equal to $0$ on $\partial D$). Under mild assumptions
on $D$, the solution to this problem exists, is unique, and
$$
u(x)=\int_{D}dy G_{D}(x,y)g(y)=E_x(\int_0^{\tau}dt f(Z_t)).
$$

Not surprisingly, the Poisson kernel is closely related to the
Green's function. More precisely, if $n=n_{z,D}$ is the inwards pointing normal vector at $z \in \partial D$,
then, as $\epsilon$ goes to $0$,
$$G_D (x, z + \epsilon n ) \sim 2 \epsilon h_D (x,z).$$

In the case of the unit disc $\U :=\{x: |x|<1\}$ in the complex plane, the
Poisson kernel and the Green's function can be explicitly computed:
$$
h_{\U}(x,z)=\frac{1-|x|^{2}}{ 2 \pi |x-z|^{2}}
$$
and
$$
G_{\U}(x,y)=\frac{-1}{\pi}\log\frac{|x-y|}{|1-x\bar{y}|}
$$
for $x\in \U, y\in \U$, and $z\in\partial \U$.

\subsection{Conformal invariance}
Conformal invariance of planar Brownian motion, first observed by Paul L\'evy \cite {Levy}, can be described as
follows: if one considers a planar Brownian motion $Z$ started from
$x$ and stopped at its first exit time of a simply connected domain
$D$, and if $\Phi$ denotes a conformal map from $D$ onto some other
domain $D'$, then the law of $\Phi(Z)$ is that of a Brownian motion
started from $\Phi(x)$ and stopped at its first exit time of $D'$.
Actually, for this statement to be fully true, one has to
reparametrize time of $\Phi(Z)$ in a proper way. The rigorous
statement is that for all $t <\tau$,
$$
\Phi(Z_{t})=Z'_{H_{t}} \hbox { with }
H_{t}=\int_{0}^{t}ds|\Phi'(Z_{s})|^{2},
$$
where $Z'$ is a Brownian motion started from $\Phi(x)$, stopped
at $\tau'= H_\tau$, which is its exit
time of $D'$.

Conformal invariance of Brownian motion is closely related to the
conformal invariance of the Green's function and of the Poisson
kernel.
Let us give a rather convoluted explanation of the conformal invariance of Green's functions using Brownian motion (a direct proof using the
analytic characterization of the Green's function is much more straightforward) that will be helpful for what follows.
 Suppose that $x$ and $y$ are in $D$ and that $\epsilon$ is
very small. We have seen that the expected time spent in the ball
$U(y,\epsilon)$, centered at $y$ and of radius $\epsilon$, by the
Brownian motion $Z$ started at $x$ behaves like
$$
\pi\epsilon^{2}G_D (x,y)
$$
when $\epsilon\to 0$.
Equivalently, the expected time spent in the ball $U(\Phi(y),|\Phi'(y)|\epsilon)$ by
 the
Brownian motion $\beta$ started at $\Phi(x)$, behaves like
$$
\pi|\Phi'(y)|^{2}\epsilon^{2} G_{D'}(\Phi(x),\Phi(y))
$$
as $\epsilon\to 0$. The process $\Phi(Z)$ can be viewed as a
time-changed Brownian motion, and the time-change when $Z$ is close
to $y$ is described via $H_{t}$. It follows easily that this expected
time of $\Phi(Z)$ spent in the ball $U(\Phi(y),|\Phi'(y)|\epsilon)$
behaves like
$$\frac{\pi|\Phi'(y)|^{2}\epsilon^{2} G_{D'}(\Phi(x),\Phi(y))}{|\Phi'(y)|^{2}}= \pi \epsilon^2 G_{D'}(\Phi(x),\Phi(y)).
$$
As a result, we have indeed that
\begin{equation}\label{conf inv Green}G_{\Phi(D)}(\Phi(x),\Phi(y))=G_{D}(x,y).
\end{equation}
For a more rigorous derivation along the same lines, we can use the integral
representation of occupation times of domains : on the one hand,
\begin{eqnarray*}
{E_{\Phi(x)}(\int_{0}^{\tau_{D'}}dtf(Z'_{t}))}
&=&\int_{D'}dy G_{D'}(\Phi(x),y)f(y)\\
&=&\int_{D}|\Phi'(y)|^{2}dy G_{D'}(\Phi(x),\Phi(y))f(\Phi(y))
\end{eqnarray*}
for indicator functions $f=1_{A}$, and on the other hand,
\begin{eqnarray*}
{E_{\Phi(x)}(\int_{0}^{\tau_{D'}}dtf(Z'_{t}))}
&=&E_{x}(\int_{0}^{\tau_{D}}|\Phi'(Z_{t})|^{2}dt
f(\Phi(Z_{t})))\\
&=&\int_{D}dy G_{D}(x,y)f(\Phi(y))|\Phi'(y)|^{2}
.\end{eqnarray*}

Conformal invariance of planar Brownian motion can also be used in a similar way to see that
\begin{equation}
|\Phi'(z)| \ h_{\Phi(D)}(\Phi(x),\Phi(z))=h_{D}(x,z)
\end{equation}
for all $x\in D, z\in\partial D$ when
$\partial D$ is smooth. Let us stress again that these conformal
invariance properties of the Green's functions and of the Poisson
kernel can be derived much more directly without any reference to
Brownian paths.

Note that $G_{\U}(0,y_{0})=- \pi^{-1} \log |y_{0}|$ for all $y_{0}\neq 0$.
The formula for $G_{\U}(x,y)$ then follows immediately, using the
M\"{o}bius transformation $\phi_{x}$ of $\U$ onto itself that maps $x$ onto $0$ and vice-versa
(this is the map $z\mapsto (z-x)/(1-\bar{x}z)$) because then $G_{\U}(x,y)=G_{\U}(0,\phi_{x}(y))$.
Note also that this conformal invariance also provides one
possible explanation of the symmetry of the Green's function $G_{\U} (x,y ) = G_{\U} (y,x)$ (because for any $x$ and $y$, there exists
a conformal map from $D$ onto itself that maps $x$ onto $y$ and $y$ onto $x$).

Similarly, since clearly $h_{\U}(0,z)=1/(2\pi)$ for all $z\in\partial
\U$, the formula for $h_{\U}(x,z)$ recalled at the end of the previous subsection follows using conformal
invariance.

\subsection{The Gaussian Free Field}

In the present text, we will briefly relate our
Brownian excursions to the Gaussian Free Field, which is a classical and basic building block in
Field theory, see for instance \cite {Nelson1973,Gawedzki}. So we recall its definition, in the
Gaussian Hilbert space framework (as in \cite {GFF2006} for instance):
Consider the space $H_s(D)$ of smooth, real-valued functions on
$\mathbb R^2$ that are supported on a compact subset of a domain $D
\subset \mathbb R^d$ (so that, in particular, their first
derivatives are in $L^2(D)$). This space can be endowed with the
\textit{Dirichlet inner product} defined by
$$
(f_1,f_2)_{\nabla} = \int_D dx(\nabla f_1 \cdot \nabla f_2)
$$
It is immediate to see that this Dirichlet inner product is invariant under conformal transformation.
Denote by $H(D)$ the Hilbert space completion of $H_s(D)$. The quantity
$(f,f)_{\nabla}$ is called the \textit{ Dirichlet energy} of $f$.

A {\em Gaussian Free Field} is any Gaussian Hilbert space
$\mathcal G(D)$ of random variables denoted by
``$(h,f)_{\nabla}$''---one variable for each $f \in H(D)$---that
inherits the Dirichlet inner product structure of $H(D)$, i.e.,
$$\mathbb E [(h,a)_{\nabla} (h, b)_{\nabla}] = (a,b)_{\nabla}.$$
In other words, the map from $f$ to the random variable $(h,
f)_{\nabla}$ is an inner product preserving map from $H(D)$ to
$\mathcal G(D)$. The reason for this notation is that it is possible
to view $h$ as a random linear operator, but we will not need this
approach. We also view $(h,\rho)$ as being well defined for all
$\rho\in(-\triangle)H(D)$ (if $\rho=-\triangle f$ for some $f\in
H(D)$, then we denote $(h,\rho)=(h,f)_{\nabla}$).

When $\rho_1$ and $\rho_2$ are in $H_s(D)$, the covariance of $(h,
\rho_1)$ and $(h, \rho_2)$ can be written as
$(-\triangle^{-1}\rho_{1},-\triangle^{-1}\rho_{2})_{\nabla}=(\rho_1,
-\Delta^{-1} \rho_2)= (-\Delta^{-1} \rho_1, \rho_2)$. From the
Poisson problem that we discussed before, $-\Delta^{-1} \rho$ can be
written using the Green's function as
$$[-\Delta^{-1}\rho](x) = \frac{1}{2}\int_{D}dy \ G_{D}(x,y) \rho(y) ,$$
we may also write:
\begin{equation}\label{cov gff}
\text{Cov}[(h, \rho_1), (h,\rho_2)] = \frac{1}{2}\int dxdy
\ G_{D}(x,y)\rho_1(x) \rho_2(y)
\end{equation}

Both the Dirichlet inner product and the Gaussian Free Field inherit naturally conformal invariance properties from
the conformal invariance of the Green's
function. The $2$-dimensional Gaussian free field (GFF) is a particular rich object, in which a number of geometric features can be detected, and that
 has been shown to play
a central role in the theory of random
surfaces and conformally invariant geometric structures, see \cite {Sheffield2011} and the references therein.

\section{Occupation times of Brownian excursions}
{\bf Brownian excursion measure.}
Let us first very briefly recall the construction of Brownian
excursion measures. For the unit disc $\U$, for each $\epsilon>0$, let
$\mu_{\epsilon}$ denote the measure of total mass $1/\epsilon$ defined as $1/\epsilon$
times the law of a Brownian motion started uniformly on the circle
of radius $(1-\epsilon)$, and stopped at its first hitting time of the
unit circle. In some appropriate topology, the measures $\mu_{\epsilon}$
converge when $\epsilon\rightarrow 0$ to an infinite measure $\mu$ on
two-dimensional paths that start and end on the unit circle. For a
general simply connected domain $D$, the excursion measure $\mu_{D}$
can either be defined as the image of $\mu$ by the conformal map
$\Phi$ that maps $\U$ onto $D$, or alternatively in an analogous way
as in the disc, by integrating over the choice of the starting point
of the excursion on $\partial D$. The fact that these two
definitions are equivalent is the conformal invariance property of
the Brownian excursion measures. See e.g. \cite{wernerleshouches} for details
and references.

Note that $\mu$ is a measure on paths $(B_t, 0 < t < \tau)$ that start and end on $\partial D$ (i.e., $B_0 \in \partial D$ and $B_\tau \in \partial D$) that are ``oriented'', i.e. $B_0$ and $B_\tau$ do a priori not play the same role. However, it turns out that the Brownian excursions are reversible i.e., that $(B_t, 0 < t < \tau )$
and $(B_{\tau -t}, 0 < t < \tau)$ are defined under the same measure (this can for instance be easily seen using the definition in the case where $D$ is the upper half-plane).

\medbreak
\noindent
{\bf Brownian excursion occupation times and the Dirichlet problem.}
Let us first make a comment on the relation between the Brownian excursion measure and the Dirichlet problem. Let $u$ be the solution to the Dirichlet problem, i.e. $\Delta u=0 $ in $\U$ and $u=f$ on $\partial\U.$ For all $z\in\partial \U$ and all positive $\epsilon$, we have that
\begin{eqnarray*}
{E_{(1-\epsilon)z}(\int_0^{\tau}dt 1_A(\gamma_t)f(\gamma_{\tau}))}
&=& E_{(1-\epsilon)z}(\int_0^{\infty}dt 1_A(\gamma_t)1_{t\leq \tau}f(\gamma_{\tau}))\\
&=& E_{(1-\epsilon)z}(\int_0^{\infty}dt 1_A(\gamma_t)1_{t\leq \tau}E(f(\gamma_{\tau})|\mathcal{F}_t))\\
&=& E_{(1-\epsilon)z}(\int_0^{\tau}dt 1_A(\gamma_t)E_{\gamma_t}(f(\gamma_{\tau})))\\
&=& E_{(1-\epsilon)z}(\int_0^{\tau}dt 1_A(\gamma_t)u(\gamma_t))\\
&=& \int_A dy G_{\U}((1-\epsilon)z, y)u(y)
\end{eqnarray*}
And for the Brownian excursion measure $\mu=\mu_{\U}$, we have that
\begin{eqnarray*}
{\mu(\int_0^{\tau}dt 1_A(\gamma_t)f(\gamma_{\tau}))}
&=& \lim_{\epsilon\rightarrow 0} \int_0^{2\pi}\frac{d\theta}{\epsilon}E_{(1-\epsilon)e^{i\theta}}(\int_0^{\tau}dt 1_A(\gamma_t)f(\gamma_{\tau}))\\
&=& \lim_{\epsilon\rightarrow 0} \int_0^{2\pi}\frac{d\theta}{\epsilon}\int_A dy G_{\U}((1-\epsilon)e^{i\theta}, y)u(y)\\
&=& \int_0^{2\pi}2d\theta\int_A dy\ h_{\U}(y,e^{i\theta})u(y)\\
&=& 2\int_A dy\ u(y)\int_0^{2\pi}d\theta h_{\U}(y,e^{i\theta})\\
&=& 2\int_A dy\ u(y)
\end{eqnarray*}
That is to say, we can represent the solution to the Dirichlet problem via the Brownian excursion measure by the formula
$$
\mu(\int_0^{\tau}dt 1_A(\gamma_t)f(\gamma_{\tau}))=2\int_A dy\ u(y)
$$
Since the Brownian excursion is reversible, we also have that
\begin{equation}\label{Dirichlet and exc}
\mu(f(\gamma_{0})\int_0^{\tau}dt 1_A(\gamma_t))=2\int_A dy\ u(y)
\end{equation}
Hence, if we put a weight $f$ on starting point of the excursion, then the mean occupation time spent in $A$ is measured by the integral of $u$ on $A$, where $u$ is the solution to the corresponding Dirichlet problem. By conformal invariance, (\ref{Dirichlet and exc}) also holds for any simply connected domain.

We would like to note that,
if we set $f=1$ in (\ref{Dirichlet and exc}), we get that $\mu_{D}(\int_0^\tau dt 1_{A} ( \gamma_t))$ is equal to twice the area of $A$. In particular, $\mu_{D} (\tau)$ is therefore just twice the area of $D$.

\medbreak
\noindent
{\bf The covariance structure.}
We now turn our attention towards the proof of ($\ref{exc cov}$). This formula can
be understood as follows: we can cut $A\times B$ into very small
pieces, calculate on each small piece and then add all these pieces
together. On each small piece $dx \times dy$, the Brownian excursion
starts from the boundary, firstly it hits the small piece $dx$
(with a small probability), after this time, it is a true Brownian
motion starting nearby $x$, which is (almost) independent of the
past and then the expected time of this new Brownian motion spent in
the neighborhood of $y$ before exiting $D$ is close to
$G_{D}(x,y)dy$. When we add up all these small pieces together
and we obtain the right-hand side of the formula.

For a precise calculation, we first consider the case where $D=\U$ as
the general case will then follow from conformal invariance. We also use
the notation that $\mu=\mu_{\U}$. Let $\gamma$ denote a
 Brownian excursion in $\U$. For all $z\in\partial \U$ and all positive $\epsilon$,
\begin{eqnarray*}
\lefteqn{E_{(1-\epsilon)z}(\int_{0}^{\tau}ds
1_{A}(\gamma_{s})\int_{s}^{\tau}dt 1_{B}(\gamma_{t}))}\\
&=&E_{(1-\epsilon)z}(\int_{0}^{\tau}ds 1_{A}(\gamma_{s})E(\int_{s}^{\tau}dt 1_{B}(\gamma_{t})|\mathcal {F}_{s}))\\
&=&E_{(1-\epsilon)z}(\int_{0}^{\tau}ds
1_{A}(\gamma_{s})E_{\gamma_{s}}(\int_{0}^{\tau}dt 1_{B}(\gamma_{t})))\\
&=&E_{(1-\epsilon)z}(\int_{0}^{\tau}ds
1_{A}(\gamma_{s})G_{\U}(\gamma_{s},B))\\
&=&\int_{A}dy G_{\U}((1-\epsilon)z,y)G_{\U}(y,B).
\end{eqnarray*}
And for the Brownian excursion measure, we have that
\begin{eqnarray*}
\lefteqn{\mu(\int_{0}^{\tau}ds 1_{A}(\gamma_{s})\int_{s}^{\tau}dt 1_{B}(\gamma_{t}))}\\
&=&\lim_{\epsilon\rightarrow 0}\int_0^{2\pi}\frac{d\theta}{\epsilon } E_{(1-\epsilon)e^{i\theta}}(\int_{0}^{\tau}ds 1_{A}(\gamma_{s})\int_{s}^{\tau}dt 1_{B}(\gamma_{t}))\\
&=&\lim_{\epsilon\rightarrow 0}\int_0^{2\pi}\frac{d\theta}{\epsilon}\int_{A}dy G_{\U}((1-\epsilon)e^{i\theta},y)G_{\U}(y,B)\\
&=&\int_0^{2\pi}2d\theta\int_{A}dy h_{\U}(y,e^{i\theta})G_{\U}(y,B)\\
&=&2\int_{A}dy G_{\U}(y,B)\int_0^{2\pi}d\theta h_{\U}(y,e^{i\theta})\\
&=&2\int_{A}dy G_{\U}(y,B)
\end{eqnarray*}

By symmetry of the Green's function ($G_{\U}(x,y)=G_{\U}(y,x)$), we
have that
$$
\mu(\int_{0}^{\tau}ds 1_{A}(\gamma_{s})\int_{0}^{\tau}ds
1_{B}(\gamma_{s}))=4\int_{A\times B}dxdy
G_{\U}(x,y).
$$
This concludes the proof of the equation (\ref{exc cov}), since we can use to conformal invariance to derive the formula for general simply connected domain $D$. More
generally, we have that
\begin{equation}
\mu_{D}(\int_{0}^{\tau}ds f(\gamma_{s})\int_{0}^{\tau}ds
g(\gamma_{s}))=4\int dx\ dy G_{D}(x,y)f(x)g(y)
\end{equation}
for all measurable bounded functions $f$ and $g$.

\medbreak
\noindent
{\bf Large intensity clouds of excursions and GFF.}
Let us now use this formula to make a link between Brownian
excursions and the GFF. For this we are going to use Poissonian
cloud of excursions in $D$, as in \cite{conf2005}. Recall that a
Poisson cloud of excursions with intensity $c\mu_{D}$ is a random
countable family of Brownian excursions in $D$, that is defined as a
Poisson point process with intensity $c\mu_{D}$.

In particular, the union of two independent Poissonian clouds of
Brownian excursions in $D$ with intensity $c_{1}\mu_{D}$ and
$c_{2}\mu_{D}$ is a Poissonian cloud of excursions in $D$ with
intensity $(c_{1}+c_{2})\mu_{D}$.

Let us now consider an i.i.d. sequence $M^{j},j\geq 1$ of Poissonian
clouds of excursions in $D$ with the common intensity $\mu_{D}$. For
each $j\geq 1$, and each $f\in (-\Delta)H(D)$, define the
``cumulative occupation'' time of $M^{j}$ by
$$
X^{j}_{f}=\sum_{\gamma\in M^{j}}\int_{0}^{\tau(\gamma)}ds f(\gamma_{s}).
$$
The fact that $\mu (\tau)$ is finite (as soon as the area of $D$ is finite) ensures that
$X^j_f$ is almost surely finite (as soon as $f$ is bounded) because its expectation is bounded.
We then define
$$
\tilde{X}^{j}_{f}=X_{f}^{j}-E(X_{f}^{j}).
$$
On an enlarged probability space, we can also define an i.i.d.
family of random variable $\epsilon_{\gamma}$ indexed by the set of
excursions in $\cup_{j}M^{j}$ such that
$P(\epsilon_{\gamma}=1)=P(\epsilon_{\gamma}=-1)=1/2$. We can then
define
$$
Y_{f}^{j}=\sum_{\gamma\in M^{j}} \epsilon_{\gamma}\int_{0}^{\tau(\gamma)}ds
f(\gamma_{s}).
$$
It is easy to see that $Y_{f}^{1},Y_{f}^{2},Y_{f}^{3}, \ldots$ are
i.i.d. centered random variables with common variance
$$
\sigma_{f}^{2}=\mu_{D}(\int_{0}^{\tau}ds
f(\gamma_{s})\int_{0}^{\tau}ds f(\gamma_{s}))=4\int dxdy
G_{D}(x,y)f(x)f(y).
$$
The same is true for
$\tilde{X}_{f}^{1},\tilde{X}_{f}^{2},\tilde{X}_{f}^{3}, \ldots$. By the
Central Limit Theorem, we have that
$$
\frac{1}{\sqrt{N}}(Y_{f}^{1}+...+Y_{f}^{N})
$$
converges in law as $N \to \infty$ to a centered Gaussian random
variable $Y_{f}$ with variance $\sigma_{f}^{2}$. The same holds for
the sequence
$(\tilde{X}_{f}^{1}+...+\tilde{X}_{f}^{N})/\sqrt{N}$.

Hence, we see that the GFF can be viewed as the limit (in law, and in
the sense of finite-dimensional distributions) of the occupation
times fluctuations of a Poisson cloud of Brownian excursions, when
the intensity tends to infinity.

\medbreak
\noindent
{\bf Higher-order ``moments''.}
We just mention that our proof can be adapted directly in order to show that for all $p \ge 2$:
\begin {eqnarray*}
\lefteqn {\mu_{D} ( \int_{(0,\tau)^p} dt_1 \ldots dt_p 1_{t_1 < \ldots < t_p} 1_{A_1} (\gamma_{t_1}) \cdots 1_{A_p} ( \gamma_{t_p} ) )}
\\
&=&
2 \int_{A_1 \times \cdots \times A_p} dx_1 \cdots dx_p G_{D}(x_1, x_2) \times \cdots \times G_D ( x_{p-1}, x_p)
\end {eqnarray*}
which gives for instance (when one sums over all possible order of visits) a formula for
$ \mu_{D}( (\int_0^\tau f(\gamma_s) ds)^p)$. We have chosen to focus on the case $p=2$ because of the above-mentioned link with Gaussian fields.

\medbreak
\noindent
{\bf Non-simply-connected domains.}
Suppose that $D$ is a finitely connected open domain in the plane. Then, by Koebe's uniformization Theorem (see \cite {Koebe}), it is possible to map it conformally onto a circular domain i.e., the unit disk $U$ punctured by a finite number of disjoint closed disks. It is very easy to generalize the definition of the Brownian excursion measure in
circular domains (adding the contributions corresponding to starting points in the neighborhood of each of the boundary disks), and to see that all our proofs go through without any real difficulty, so that all our statements are in fact valid also in circular domains. One can then {\em define} the excursion measure in $D$ via conformal invariance starting from circular domains, and then, by conformal invariance of all the quantities involved, we easily see that all our statements are also valid in $D$.

\section{Occupation times of Brownian loops}
\noindent
{\bf Brownian loop measure.}
We now briefly recall the construction of
the Brownian loop measure \cite{soup2004}. As for the Brownian excursion measure, we can first define it in
the unit disc, and then define it in any other simply connected domain using conformal invariance (and one then checks that this is indeed consistent with
other possible constructions).

For any $r \in (0,1]$, define $U_r = r\U$. For any $x \in U_r$ and any $z \in \partial U_r$, one can define the
Brownian motion started at $x$ and conditionned to exit $U_r$ at $z$ (this can be rigorously defined as the limit when $\epsilon \to 0$
of the law of the Brownian motion conditioned to exit $U_r$ in an $\epsilon$-neighborhood of $z$). Let us denote this probability measure
by $P_{x \to z}^r$.
Then, as for the excursion measure, one can let $x \to z$, and renormalize it in order to get a measure on macroscopic sets i.e. define
$$
m_z^r (\cdot)  = \lim_{\epsilon \to 0} \epsilon^{-1} h_{U_r}(z+\epsilon n,z)P_{z+\epsilon n \to z}^r ( \cdot)
$$
where $n=n_{z,U_r}$ is the inwards pointing normal vector at $z\in \partial U_r.$
Then, one can define the loop measure in $\U$ by integrating $z$ on $\partial U_r$, and then integrating $r$ from $0$ to $1$:
$$
\lambda_{\U} (\cdot) =\int_{0}^{1}r dr\int_{0}^{2\pi}d\theta\  m^r_{re^{i\theta}}(\cdot).
$$

In fact, the above definition is not quite the loop measure because it defines a measure on parametrized loops. We will forget about the precise parametrization of the loop and view $\lambda_{\U}$ as a measure on loops defined modulo monotone reparametrization (where the time-parameter should be viewed as an element of the circle, because the end-point of the loop is the same as the starting point, this is possible). It turns out that this definition of $\lambda_{\U}$ is then invariant under the Moebius transformations that map the unit disc onto itself. Hence, it is possible to define,
for a general simply connected domain $D$, the loop measure $\lambda_{D}$ as the image of $\lambda_{\U}$ by any conformal map $\Phi$
that maps $\U$ onto $D$. And we usually denote $\lambda=\lambda_{\U}.$

Before going on, we would like to say a word on the value of $\lambda(\tau)$. In fact, by direct computation we have that $\lambda(\tau)=\infty$ which is very different from $\mu(\tau)$ mentioned before. A direct way to check that $\lambda( \tau) = \infty$ goes as follows. Consider $D$ to be the square $[0,1]^2$. For any dyadic square $d$ in $D$ with sidelength $2^{-n}$, a direct scaling argument shows that the mass (for $\lambda$) of the set of loops that stay in $d$ and have a time-length in $[4^{-n} ,2 \times 4^{-n})$  does not depend on $d$. Hence, if we sum this quantity over all dyadic squares $d$ in $D$, and because
$\sum_n 4^n 4^{-n} = \infty$, we readily see that $\lambda (\tau) = \infty$.

However, almost the same argument ensures that $\lambda ( \tau^{1+\epsilon})$ is finite for $\epsilon > 0$ (and bounded $D$). Indeed, in the case of the unit square, we can decompose the set of loops with time-length in $[4^{-n}, 4^{1-n})$ according to the dyadic square in which its lowest point lies. This leads readily to the bound
$$ \lambda (1_{\tau < 1 }  \tau^{1+ \epsilon}) \le C \sum_{n \ge 1} 4^n (4^{1-n})^{1+ \epsilon} < \infty$$
and one can see by other means that $\lambda ( \tau > t )$ decays exponentially fast as $t \to \infty$.
In particular, we get that $\lambda ( \tau^2 )$ is finite (as soon as $D$ is bounded).

\medbreak
\noindent
{\bf Covariance structure.}
Our goal is now to prove (\ref {loop cov}). As before, we are going to derive the result first in the case where $D=\U$, and the general result will then follow using conformal invariance.
Again, it will be convenient to (loosely speaking) divide $A\times B$ into infinitesimal pieces $dx\times dy$, make the
computation on each piece, and then add all these pieces together. Clearly, this will
give a formula of the type
$$
\lambda_{D}(\int_{0}^{\tau}ds 1_{A}(\gamma_{s})\int_{0}^{\tau}ds
1_{B}(\gamma_{s}))=\int_{A\times B}dxdy F_{D}(x,y)
$$
where $F_{D}(x,y)$ is the ``covariance'' function between $x$ and $y$
determined by the Brownian loop measure. Just as what we have done
to derive the conformal invariance of the Green's function in the
equation (\ref{conf inv Green}), we can also derive the conformal
invariance of $F$:
$$
F_{\Phi(D)}(\Phi(x),\Phi(y))=F_{D}(x,y).
$$

To determine $F_{D}(x,y)$, it is enough
to describe $F_{\U}(0,y_{0})$ for  $y_{0}\in(0,1)$, because there exists
a $y_0$ and a conformal map $\Phi$ from $D$ onto $\U$ such that
$\Phi(x)=0,\Phi(y)=y_{0}$.

And now begin our computation. For $r\in(0,1), z\in\partial U_r$, we can write
\begin{eqnarray*}
\lefteqn{E^r_{z+\epsilon n\to z}(\int_0^{\tau}ds 1_A(\gamma_s)\int_s^{\tau}dt 1_B(\gamma_t))}\\
&=& \lim_{\epsilon'\to 0} \Bigl( \frac{1}{P^r_{z+\epsilon n}(\gamma_{\tau}\in U(z,\epsilon')\cap\partial U_r)} \\
&& \quad \times  E^r_{z+\epsilon n}(\int_0^{\tau}ds 1_A(\gamma_s)\int_s^{\tau}dt 1_B(\gamma_t)1_{\gamma_{\tau}\in U(z,\epsilon')\cap\partial U_r}) \Bigr) \\
&=& \lim_{\epsilon'\to 0} \Bigl( \frac{1}{h_{U_r}(z+\epsilon n,U(z,\epsilon')\cap\partial U_r)} \\
&& \quad \times E^r_{z+\epsilon n}(\int_0^{\tau}ds 1_A(\gamma_s)\int_s^{\tau}dt 1_B(\gamma_t)h_{U_r}(\gamma_t,U(z,\epsilon')\cap\partial U_r)) \Bigr) \\
&=& \frac{1}{h_{U_r}(z+\epsilon n,z)}E^r_{z+\epsilon n}(\int_0^{\tau}ds 1_A(\gamma_s)\int_s^{\tau}dt 1_B(\gamma_t)h_{U_r}(\gamma_t,z)).
\end{eqnarray*}
Hence,
\begin{eqnarray*}
\lefteqn{h_{U_r}(z+\epsilon n, z)E^r_{z+\epsilon n\to z}(\int_0^{\tau}ds 1_A(\gamma_s)\int_s^{\tau}dt 1_B(\gamma_t))}\\
&=& E^r_{z+\epsilon n}(\int_0^{\tau}ds 1_A(\gamma_s)\int_s^{\tau}dt 1_B(\gamma_t)h_{U_r}(\gamma_t,z))\\
&=& \int_A dx\ G_{U_r}(z+\epsilon n, x)\int_B dy\ G_{U_r}(x,y)h_{U_r}(y,z)
\end{eqnarray*}
and letting $\epsilon \to 0$, we get
\begin{eqnarray*}
\lefteqn{m^r_z(\int_0^{\tau}ds 1_A(\gamma_s)\int_0^{\tau}ds 1_B(\gamma_s))}\\
&=& \lim_{\epsilon\to 0}\frac{2}{\epsilon}\int_A dx\ G_{U_r}(z+\epsilon n, x)\int_B dy\ G_{U_r}(x,y)h_{U_r}(y,z)\\
&=& 4\int_{A\times B} dxdy G_{U_r}(x,y)h_{U_r}(x,z)h_{U_r}(y,z).
\end{eqnarray*}
For simplicity, we define a new kernel
$$
K_{U_r}(x,y)=4\int_0^{2\pi}d\theta h_{U_r}(x,re^{i\theta})h_{U_r}(y,re^{i\theta})
$$
and then we have that
$$
\lambda_{\U}(\int_0^{\tau}ds 1_A(\gamma_s)\int_0^{\tau}ds 1_B(\gamma_s))=\int_0^1 rdr \int_{A\times B}dxdy G_{U_r}(x,y)K_{U_r}(x,y).
$$
Note that
$$
K_{U_r}(rx,ry)=\frac{1}{r^2}K_{\U}(x,y)
$$
and
$$
G_{U_r}(rx,ry)=G_{\U}(x,y).$$
Furthermore,
$ K_{\U}(0,y)={2/\pi}$.

Suppose that  $A= U(0,\epsilon)$ and $B=
U(y_{0},\delta)$ where $\epsilon$ and $\delta$  are both small.
In our decomposition of $\lambda_{\U}$, the loop can visit $B$ only if it started on a circle of radius $r > y_0$.
Hence,
on the one hand, as $\epsilon$ and $\delta$ tend to $0$,
\begin{eqnarray*}
\lefteqn{\lambda_{\U}(\int_{0}^{\tau}ds
1_{A}(\gamma_{s})\int_{0}^{\tau}ds
1_{B}(\gamma_{s}))}\\
&=&\int_{y_0}^{1}rdr\int_{(A\bigcap U_{r})\times (B\bigcap
U_{r})}dxdy
G_{U_{r}}(x,y)K_{U_{r}}(x,y)\\
&\sim& \int_{y_0}^{1}rdr
(\pi\epsilon^{2}\pi\delta^{2})G_{U_{r}}(0,y_{0})K_{U_{r}}(0,y_{0})\\
&=&(\pi\epsilon^{2}\pi\delta^{2})\int_{y_{0}}^{1}\frac{1}{r}dr
G_{\U}(0,\frac{y_{0}}{r})K_{\U}(0,\frac{y_{0}}{r})\\
&=&(\pi\epsilon^{2}\pi\delta^{2})\frac{2}{\pi^{2}}\int_{y_{0}}^{1}dr(-\frac{1}{r}\log(\frac{y_{0}}{r}))\\
&=&(\pi\epsilon^{2}\pi\delta^{2})\frac{1}{\pi^{2}}(\log y_{0})^{2}
\end{eqnarray*}

On the other hand, this quantity is precisely behaving as $(\pi\epsilon^{2}\pi\delta^{2})F_{\U}(0,y_{0})$ and
as a result, we get that
$$
F_{\U}(0,y_{0})=\frac{1}{\pi^{2}}(\log
y_{0})^{2}=(G_{\U}(0,y_{0}))^2.
$$
We can then conclude that (\ref {loop cov}) holds in $\U$, and then also in $D$ by conformal invariance.
More generally, we have that
$$
\lambda_{D}(\int_{0}^{\tau}ds f(\gamma_{s})\int_{0}^{\tau}ds
g(\gamma_{s}))=\int_{A\times B}dxdy
(G_{D}(x,y))^2f(x)g(y).
$$
for all measurable bounded functions $f$ and $g$.

\medbreak
\noindent
{\bf Brownian loop-soups and fields.}
Just as in the case of Brownian excursion measure, we can use this
formula to make a link between Brownian loops and some Gaussian Fields.
Let $M^{j},j\geq 1$ be a sequence of i.i.d Poissonian clouds of
loops in $D$ with the common intensity $\lambda_{D}$. We can try to give
the same definitions of the quantities $\tilde{X}_{f}^{j},
Y_{f}^{j}, j\geq 1$. However, things are a little more complicated, due to the fact that
the same scaling argument that showed that $\lambda ( \tau) = \infty$ implies that
$$ \sum_{\gamma \in M^j} \tau(\gamma) = \infty $$
almost surely, so that some care is needed.

The definition of $Y_f^j$ is however not a big problem. Recall that on an enlarged probability space, one associates to each loop $\gamma$ a
random variable $\eps_\gamma$ with $E( \eps_\gamma) = 0$ and $E ( (\eps_\gamma)^2 ) = 1$.
But equation (\ref {loop cov}) precisely ensures that the sum
$$ \sum_{\gamma \in M^j} \eps_\gamma \int_0^{\tau(\gamma)} f ( \gamma_s) ds $$
makes sense in $L^2$, and that its second moment is equal to
$$
\sigma^{2}_{f}=\lambda_{D}(\int_{0}^{\tau}dsf(\gamma_{s})\int_{0}^{\tau}dsf(\gamma_{s}))=\int
dx dy (G_{D}(x,y))^2f(x)f(y)
$$
which is finite.

Then, just as in the case of the clouds of excursions, the sequence $Y_{f}^{1}$,
$Y_{f}^{2}, \ldots$ is made of i.i.d centered random variables with
common variance $\sigma_f^2$.
By the Central Limit Theorem,
$$
\frac{1}{\sqrt{N}}(Y_{f}^{1}+...+Y_{f}^{N})
$$
converges in law as $N \to \infty$ to a centered Gaussian random
variable with variance $\sigma_{f}^{2}$. Hence, we obtain
another Gaussian Field, characterized by this new covariance structure.

It is also still possible to make sense of $\tilde X_f^j$ even though it is not possible to define $X_f^j$. It suffices to partition the set of loops (in $D$) into a countable set of loops $A_k, k \ge 1$ such that for each $k$, $\lambda ( \tau 1_{\gamma \in A_k})$ is finite (for instance, one can take
$A_k = \{ \gamma \ : \ \tau ( \gamma) > 1/k \} \setminus (A_1 \cup \ldots \cup A_{k-1}) \}$. Then, one can define
$$
\tilde X^j_f = \sum_{k \ge 1}  \left( \sum_{\gamma \in A_k \cap M^j} \int_0^{\tau (\gamma)} f( \gamma_s) ds
- E ( \sum_{\gamma \in A_k \cap M^j} \int_0^{\tau (\gamma)} f( \gamma_s) ds) \right)
$$
and check that this sum with respect to $k$ converges in $L^2$, and that its second moment is the same as that of $Y_f$.
The rest of the argument is again the same.

\section{Intersections of Brownian excursions}

In this section, we try to find the relation between intersections of Brownian
excursion ``occupations times'' and Brownian loop occupation times, the former being defined via the intersection local time.

Let us first recall some features of Brownian intersection local times. Let $p\geq 2$
be an integer, and let $Z^{1},..., Z^{p}$ denote $p$ independent
Brownian motions in $\mathbb{R}^{2}$, started at $x^{1},..., x^{p}$
respectively. The intersection local time of $Z^{1},..., Z^{p}$ is a
random measure $\alpha(ds_{1}...ds_{p})$ on $\mathbb{R}_{+}^{p}$,
supported on
\begin{displaymath}
\{(s_{1},..., s_{p})\in\mathbb{R}_{+}^{p}:
Z_{s_{1}}^{1}=...=Z_{s_{p}}^{p}\}.
\end{displaymath}

The basic description concerning the intersection local time that we
will use goes as follows (see \cite {LeGallsaintflour} for details):
\begin{proposition} Almost surely, one can define a (random) measure $\alpha(ds_{1}...ds_{p})$ on
$\mathbb{R}_{+}^{p}$ such that, for any $A^{1},...,A^{p}$ bounded
Borel subsets of $\mathbb{R}_{+}$,
$$
\alpha(A^{1}\times...\times A^{p})=\lim_{\epsilon\to 0}
\alpha_{\epsilon}(A^{1}\times...\times A^{p})
$$
in the $L^{n}-$norm, for any $n<\infty$, where
$$
\alpha_{\epsilon}(ds_{1}...ds_{p})=ds_{1}...ds_{p}\int_{\mathbb{R}^{2}}
dy
\delta_{y}^{\epsilon}(Z^{1}_{s_{1}})...\delta_{y}^{\epsilon}(Z^{p}_{s_{p}})
$$
with
$\delta_{y}^{\epsilon}(z)=\frac{1}{\pi\epsilon^{2}}1_{U(y,\epsilon)}(z)$.
\end{proposition}

Let us use this in the context of the Brownian excursion measure. This time
we shall consider two Brownian excursions $\gamma$ and $\gamma'$
defined under the (infinite) measure $\mu_{D}\otimes\mu_{D}$, and study the
behavior of their intersection local time that spent in two disjoint sets $A$ and $B$, as before:
\begin{eqnarray*}
\lefteqn{\mu_{D}\otimes\mu_{D}(\int_{0}^{\tau}\int_{0}^{\tau'}\alpha(dtdt')
1_{(\gamma_{t}=\gamma'_{t'}\in
A)}\int_{0}^{\tau}\int_{0}^{\tau'}\alpha(dsds')
1_{(\gamma_{s}=\gamma'_{s'}\in B)})}\\
&=&\lim_{\epsilon\to
0}\mu_{D}\otimes\mu_{D}(\int_{0}^{\tau}\int_{0}^{\tau'}\alpha_{\epsilon}(dtdt')
1_{(\gamma_{t}\in A)}1_{(\gamma'_{t'}\in
A)} \\
&&\qquad\qquad\qquad
\int_{0}^{\tau}\int_{0}^{\tau'}\alpha_{\epsilon}(dsds')
1_{(\gamma_{s}\in B)}1_{(\gamma'_{s'}\in
B)})\\
&=&\lim_{\epsilon\to
0}\mu_{D}\otimes\mu_{D}(\int_{0}^{\tau}\int_{0}^{\tau'}dtdt'\int
dx
\delta_{x}^{\epsilon}(\gamma_{t})\delta_{x}^{\epsilon}(\gamma'_{t'})
1_{(\gamma_{t}\in A)}1_{(\gamma'_{t'}\in
A)}\\
&&\qquad\qquad\qquad\int_{0}^{\tau}\int_{0}^{\tau'}dsds'\int dy
\delta_{y}^{\epsilon}(\gamma_{s})\delta_{y}^{\epsilon}(\gamma'_{s'})
1_{(\gamma_{s}\in B)}1_{(\gamma'_{s'}\in
B)})\\
&=&\lim_{\epsilon\to 0} \int dx\int dy \quad
\mu_{D}\otimes\mu_{D} (\int_{0}^{\tau}dt
\delta_{x}^{\epsilon}(\gamma_{t})1_{(\gamma_{t}\in
A)}\int_{0}^{\tau}ds\delta_{y}^{\epsilon}(\gamma_{s})1_{(\gamma_{s}\in
B)}\\
&&\qquad\qquad\qquad\qquad\qquad\qquad \int_{0}^{\tau'}dt'
\delta_{x}^{\epsilon}(\gamma'_{t'})1_{(\gamma'_{t'}\in
A)}\int_{0}^{\tau'}ds'\delta_{y}^{\epsilon}(\gamma'_{s'})1_{(\gamma'_{s'}\in
B)})\\
&=&\lim_{\epsilon\to 0} \int dx\int dy
(4\int_{A\times B} dadb
\delta_{x}^{\epsilon}(a)\delta_{y}^{\epsilon}(b)G_{D}(a,b))^{2}\\
&=&16 \int _{A\times B} dxdy  (G_{D}(x,y))^2
\end{eqnarray*}

Hence, we see that pairs of Brownian excursions give rise to the same covariance structure as the Brownian loops. In a way, this is not too surprising, as for two points $x$ and $y$ that are both visited by $\gamma$ and by $\gamma'$, one sees in a way a loop structure (the part of $\gamma$ from $x$ to $y$, and then the part of $\gamma'$ back from $y$ to $x$).

Note that by a similar calculation, one gets that for any $p\geq 3$, if one defines for any $A$,
$$ T_p (A; \gamma^1, \ldots , \gamma^p) = \int_{0}^{\tau_{1}} \ldots \int_{0}^{\tau_{p}} \alpha(dt_{1} \ldots dt_{p})1_{(\gamma^{1}_{t_{1}}= \cdots =\gamma^{p}_{t_{p}}\in
A)},$$
then
$$
 \mu_{D}^{\otimes
p} ( T_p (A) T_p (B) ) =  4^{p}\int_{A\times B}dxdy(G_{D}(x,y))^p
.$$

\end{document}